\documentclass[a4paper]{amsart}

\RequirePackage{amsmath}
\RequirePackage{bm}
\RequirePackage{amssymb}
\RequirePackage{upref}
\RequirePackage{amsthm}
\RequirePackage{enumerate}
\RequirePackage{pb-diagram}
\RequirePackage{amsfonts}
\RequirePackage[mathscr]{eucal}
\RequirePackage{verbatim}
\RequirePackage{xr}
\RequirePackage{graphicx}
\usepackage{calc}
\usepackage{xspace}
\RequirePackage{color}
\RequirePackage{ifthen}

\newcommand{\cf}{cf.\@\xspace}
\newcommand{\resp}{resp.\@\xspace}

\newcommand{\al}{\alpha}
\newcommand{\bet}{\beta}
\newcommand{\ga}{\gamma}
\newcommand{\de}{\delta }
\newcommand{\e}{\epsilon}

\newcommand{\f}{\varphi}
\newcommand{\h}{\eta}

\newcommand{\ka}{\kappa}
\newcommand{\lam}{\lambda}

\newcommand{\n}{\nu}

\newcommand{\vt}{\vartheta}

\newcommand{\s}{\sigma}
\newcommand{\x}{\xi}

\newcommand{\C}{\varGamma}

\newcommand{\F}{\varPhi}

\newcommand{\Om}{\varOmega}

\newcommand{\di}[1]{#1\nobreakdash-\hspace{0pt}dimensional}

\newcommand{\fv}[2]{#1\hspace{0pt}_{|_{#2}}}

\newcommand{\so}{{\mc S_0}}

\newcommand{\const}{\tup{const}}

\newcommand{\msp[1]}[1]{\mspace{#1mu}}

\newcommand{\R}[1][n+1]{{\protect\mathbb R}^{#1}}

\newcommand{\Ss}[1][n+1]{{\protect\mathbb S}^{#1}}

\newcommand{\N}{{\protect\mathbb N}}
\newcommand{\Q}{{\protect\mathbb Q}}

\newcommand{\eR}{\stackrel{\lower1ex \hbox{\rule{6.5pt}{0.5pt}}}{\msp[3]\R[]}}
\newcommand{\eN}{\stackrel{\lower1ex \hbox{\rule{6.5pt}{0.5pt}}}{\msp[1]\N}}
\newcommand{\eO}{\stackrel{\lower1ex \hbox{\rule{6pt}{0.5pt}}}{\msc O}}

\DeclareMathOperator{\graph}{graph}

\DeclareMathOperator{\osc}{osc}

\newcommand\ra{\rightarrow}

\newcommand\pa{\partial}
\newcommand\pde[2]{\frac {\partial#1}{\partial#2}}
\newcommand\pdc[3]{\frac {\partial#1}{\partial#2_#3}}   
\newcommand\pdm[4]{\frac {\partial#1}{\partial#2_#3^#4}}   
 
\newcommand\df[2]{\frac {d#1}{d#2}}

\newcommand{\un}{\infty}
\newcommand{\A}{\forall}

\newcommand{\set}[2]{\{\,#1\colon #2\,\}}
\newcommand{\uu}{\cup}
\newcommand{\ii}{\cap}
\newcommand{\uuu}{\bigcup}

\newcommand{\uud}{ \stackrel{\lower 1ex \hbox {.}}{\uu}}
\newcommand{\uuud}[1]{ \stackrel{\lower 1ex \hbox {.}}{\uuu_{#1}}}
\newcommand\su{\subset}

\newcommand{\sminus}[1][28]{\raise 0.#1ex\hbox{$\scriptstyle\setminus$}}

\newcommand{\wed}{\wedge}

\newcommand{\abs}[1]{\lvert#1\rvert}

\newcommand{\norm}[1]{\lVert#1\rVert}

\newcommand{\spd}[2]{\protect\langle #1,#2\protect\rangle}

\newcommand\ch[3]{\varGamma_{#1#2}^#3}
\newcommand\cha[3]{{\bar\varGamma}_{#1#2}^#3}

\newcommand{\riem}[4]{R_{#1#2#3#4}}
\newcommand{\riema}[4]{{\bar R}_{#1#2#3#4}}

\newcommand{\tit}{\textit}

\newcommand{\tup}{\textup}

\newcommand{\mc}{\protect\mathcal}
\newcommand{\msc}{\protect\mathscr}

\providecommand{\bysame}{\makebox[3em]{\hrulefill}\thinspace}

\newcommand{\cq}[1]{\glqq{#1}\grqq\,}

\newcommand{\bt}{\begin{thm}}
\newcommand{\bl}{\begin{lem}}
\newcommand{\bc}{\begin{cor}}
\newcommand{\bd}{\begin{definition}}
\newcommand{\bpp}{\begin{prop}}
\newcommand{\br}{\begin{rem}}
\newcommand{\bn}{\begin{note}}
\newcommand{\be}{\begin{ex}}
\newcommand{\bes}{\begin{exs}}
\newcommand{\bb}{\begin{example}}
\newcommand{\bbs}{\begin{examples}}
\newcommand{\ba}{\begin{axiom}}
\newcommand{\bas}{\begin{assumption}}

\newcommand{\et}{\end{thm}}
\newcommand{\el}{\end{lem}}
\newcommand{\ec}{\end{cor}}
\newcommand{\ed}{\end{definition}}
\newcommand{\epp}{\end{prop}}
\newcommand{\er}{\end{rem}}
\newcommand{\en}{\end{note}}
\newcommand{\ee}{\end{ex}}
\newcommand{\ees}{\end{exs}}
\newcommand{\eb}{\end{example}}
\newcommand{\ebs}{\end{examples}}
\newcommand{\ea}{\end{axiom}}
\newcommand{\eas}{\end{assumption}}

\newcommand{\bp}{\begin{proof}}
\newcommand{\ep}{\end{proof}}
\newcommand{\eps}{\renewcommand{\qed}{}\end{proof}}

\newcommand{\bal}{\begin{align}}

\newcommand{\bi}[1][1.]{\begin{enumerate}[\upshape #1]}
\newcommand{\bia}[1][(1)]{\begin{enumerate}[\upshape #1]}
\newcommand{\bin}[1][1]{\begin{enumerate}[\upshape\bfseries #1]}
\newcommand{\bir}[1][(i)]{\begin{enumerate}[\upshape #1]}
\newcommand{\bic}[1][(i)]{\begin{enumerate}[\upshape\hspace{2\cma}#1]}
\newcommand{\bis}[2][1.]{\begin{enumerate}[\upshape\hspace{#2\parindent}#1]}
\newcommand{\ei}{\end{enumerate}}

\newcommand\ndots{\raise 0.47ex \hbox {,}\hskip0.06em\cdots %
     \raise 0.47ex \hbox {,}\hskip0.06em} 

\newcommand{\q}{\quad}
\newcommand{\qq}{\qquad}

\newcommand{\hp}{\hphantom}

\newcommand\nd{\noindent}

\newskip\Csmallskipamount                                                
\Csmallskipamount=\smallskipamount
\newskip\Cmedskipamount
\Cmedskipamount=\medskipamount
\newskip\Cbigskipamount
\Cbigskipamount=\bigskipamount

\newcommand\cvs{\vspace\Csmallskipamount}   
\newcommand\cvm{\vspace\Cmedskipamount}

\newskip\csa
\csa=\smallskipamount

\newskip\cma
\cma=\medskipamount

\newskip\cba
\cba=\bigskipamount

\newdimen\spt
\newcommand\citem{\cvs\advance\itemno by
1{(\romannumeral\the\itemno})\hskip3pt}
\newcommand{\bitem}{\cvm\nd\advance\itemno by
1{\bf\the\itemno}\hspace{\cma}}

\newcount\itemno
\newcommand{\lae}[1]{\label{E:#1}}
\newcommand{\lat}[1]{\label{T:#1}}
\newcommand{\lal}[1]{\label{L:#1}}

\newcommand{\lar}[1]{\label{R:#1}}

\newcommand{\lapp}[1]{\label{Proof:#1}}

\newcommand{\rt}[1]{Theorem~\ref{T:#1}}
\newcommand{\rl}[1]{Lemma~\ref{L:#1}}

\newcommand{\re}[1]{\eqref{E:#1}}
\newcommand{\frt}[1]{Theorem~\ref{T:#1} on page~\tup{\pageref{T:#1}}}
\newcommand{\frl}[1]{Lemma~\ref{L:#1} on page~\tup{\pageref{L:#1}}}
\newcommand{\frr}[1]{Remark~\ref{R:#1} on page~\tup{\pageref{R:#1}}}

\newcommand{\fre}[1]{\eqref{E:#1} on page~\tup{\pageref{E:#1}}}

\newcommand{\rpp}[1]{\pageref{Proof:#1}}

\newskip\thmskip
\thmskip=\parindent

\newskip\hsk
\newenvironment{hinw}{\labelsep=0pt\begin{list}{}{\labelsep=0pt\itemindent=0pt\labelwidth=0pt\leftmargin=\parindent\rightmargin=0pt\partopsep=\cba}%
\item\it\nopagebreak\nopagebreak}%
{\end{list}}

\newcommand\bh{\begin{hinw}}
\newcommand{\eh}{\end{hinw}}

\newtheoremstyle{normal}
  {\cba}
  {\cba}
  {}
  {\thmskip}
  {\bfseries}
  {.}
  {\hsk}
  {}

\newtheoremstyle{abschnitt}
  {\cba}
  {\cba}
  {}
  {\thmskip}
  {\bfseries}
  {.}
  {\hsk}
  {}

\newtheoremstyle{italic}
  {\cba}
  {\cba}
  {\itshape}
  {\thmskip}
  {\bfseries}
  {.}
  {\hsk}
  {}

\newtheoremstyle{aufgaben}
  {\cba}
  {\cba}
  {}
  {}
  {\normalsize\bfseries}
  {.}
  {\hsk}
  {}

\newtheoremstyle{break}
  {\cba}
  {\cba}
  {\itshape}
  {}
  {\bfseries}
  {.}
  {\newline}
  {}

\swapnumbers
\theoremstyle{italic}
\newtheorem{thm}[subsection]{Theorem}
\newtheorem{lem}[subsection]{Lemma}
\newtheorem{prop}[subsection]{Proposition}
\newtheorem{cor}[subsection]{Corollary}

\theoremstyle{normal}
\newtheorem{rem}[subsection]{Remark}
\newtheorem{definition}[subsection]{Definition}
\newtheorem{example}[subsection]{Example}
\newtheorem{examples}[subsection]{Examples}
\newtheorem{ex}[subsection]{Exercise}
\newtheorem{note}[subsection]{}
\newtheorem{axiom}[subsection]{Axiom}
\newtheorem{assumption}[subsection]{Assumption}

\theoremstyle{aufgaben}
\newtheorem{exs}[subsection]{Exercises}

\swapnumbers

\numberwithin{equation}{section}
\numberwithin{figure}{section}

\newenvironment{textequation}[1][0.8]
{\begin{equation}
\begin{aligned}
\begin{minipage}{#1\linewidth}}
{\end{minipage}
\end{aligned}
\end{equation}
\ignorespacesafterend}

\newcommand{\btext}{\begin{textequation}}
\newcommand{\etext}{\end{textequation}}

\def\hinweis{\@startsection{subsection}{2}%
 \z@{0.7\linespacing\@plus 0.5\linespacing}{0.7\linespacing}%
{\normalfont\itshape\indent}}

\newcounter{hours}\newcounter{minutes}
\newcommand{\printtime}{%
\setcounter{hours}{\time/60}%
\setcounter{minutes}{\time-\value{hours}*60}%
\ifthenelse{\value{minutes}<10}{\thehours :0\theminutes}{\thehours:\theminutes}}

\usepackage[german,english]{babel}
\usepackage{graphicx}
\RequirePackage{amsmath}
\RequirePackage{bm}
\RequirePackage{amssymb}
\RequirePackage{upref}
\RequirePackage{amsthm}
\RequirePackage{enumerate}
\RequirePackage{pb-diagram}
\RequirePackage{amsfonts}
\RequirePackage[mathscr]{eucal}
\RequirePackage{verbatim}
\RequirePackage{xr}
\RequirePackage{graphicx}
\usepackage{calc}
\usepackage{xspace}

\makeatletter
\RequirePackage{color}
\newcommand{\ann}[1]{\renewcommand{\@makefnmark}{\mbox{$^{\color{red}{\@thefnmark}}$}}%
\footnote {#1}}
\makeatother

\RequirePackage{upref}
\RequirePackage{amsthm}
\RequirePackage{enumerate}
\usepackage[mathscr]{eucal}

\usepackage{xr-hyper}

\usepackage{calc}

\newlength{\oddsidemarginlength}
\newlength{\topmarginlength}

\newcounter{numberoflines}
\newcounter{tempcc}
\oddsidemargin=\oddsidemarginlength
\evensidemargin=\oddsidemargin
\topmargin=\topmarginlength
\usepackage[colorlinks=true,linkcolor=blue,citecolor=blue,urlcolor=blue]{hyperref}  

\begin{document}

\flushbottom


\title[Non-scale-invariant inverse curvature flows]{Non-scale-invariant inverse curvature flows in Euclidean space}

\author{Claus Gerhardt}
\address{Ruprecht-Karls-Universit\"at, Institut f\"ur Angewandte Mathematik,
Im Neuenheimer Feld 294, 69120 Heidelberg, Germany}
\email{\href{mailto:gerhardt@math.uni-heidelberg.de}{gerhardt@math.uni-heidelberg.de}}
\urladdr{\href{http://www.math.uni-heidelberg.de/studinfo/gerhardt/}{http://www.math.uni-heidelberg.de/studinfo/gerhardt/}}
\thanks{This work has been supported by the DFG}

%
\date{\today}
%


\begin{abstract} 
We consider the inverse curvature flows $\dot x=F^{-p}\nu$  of closed star-shaped hypersurfaces in Euclidean space in case $0<p\not=1$ and prove that the flow exists for all time and converges to infinity, if $0<p<1$, while in case $p>1$, the flow blows up in finite time, and where we assume the initial hypersurface to be strictly convex. In both cases the properly rescaled flows converge to the unit sphere.
\end{abstract}

\maketitle

\tableofcontents

\setcounter{section}{0}
\section{Introduction}
We consider expanding curvature flows of star-shaped closed hypersurfaces in $\R$; these flows are also called inverse curvature flows. In \cite{cg90} we considered flows of the form
\begin{equation}\lae{1.1}
\dot x=F^{-1}\nu,
\end{equation}
where $F$ is a curvature function homogeneous of degree $1$, and proved that the flow exists for all time and converges to infinity. After a proper rescaling, the rescaled flow will converge to a sphere.

The equation \re{1.1} has the property that it is scale-invariant, i.e., if the initial hypersurface $M_0$ is star-shaped with respect to the origin and if it is scaled by a factor $\lam>0$, then the solution $x_\lam$ of equation \re{1.1} with initial hypersurface $\lam M_0$ is given by $\lam x$, where $x=x(t,\xi)$ is the original solution.

This scale-invariance seems to be the underlying reason why expanding curvature flows in Euclidean space do not develop singularities contrary to contracting curvature flows which will contract to a point in finite time, see \cite{huisken:jdg1}.

For non-scale-invariant flows, i.e., for flows satisfying
\begin{equation}\lae{1.2}
\dot x=F^{-p}\nu,\q 0<p\not=1,
\end{equation}
singularities will develop, if $p>1$. When the initial hypersurface is a sphere equation \re{1.2} is equivalent to an ODE, since the leaves of the flow will then be spheres too, and the spherical flow will develop a blow up in finite time if $p>1$. For $0<p<1$ the spherical flow will exist for all time and converge to infinity.

Flows with $0<p<1$ have been considered by Urbas \cite[Theorem 1.3]{urbas:convex} for convex hypersurfaces and he proved that the flow exists for all time, converges to infinity and that the rescaled hypersurfaces converge to a sphere. Though there are two mistakes in the proof and one in the formulation of the theorem, namely, the scale factor should have a negative exponent in equation (4.4) and in the formulation of the theorem, and  the time derivative of $\tilde H$ in equation (4.5) has the wrong sign; however,  the second mistake introduces the correct sign into the equation  which would have been resulted by choosing the correct exponent in the scale factor and by differentiating correctly, and, hence,  the proof is essentially correct.

In the case $p>1$ there are only a few special results in dimension $n=2$ and essentially only for the Gaussian curvature. Let
\begin{equation}
F=\s_n=H_n^\frac1n,
\end{equation}
where $H_n$ is the Gaussian curvature, then Schn\"urer \cite{oliver:2004surfaces} considered the case $p=2=n$ and Li \cite{li:inverse} the cases $1\le p\le 2=n$. The restriction to two dimensions is due to the method of proof for the crucial curvature estimates which relies on the monotonicity of a certain rather artificial expression depending on the principal curvatures of the flow hypersurfaces.

In this paper we consider the flow \re{1.2} for star-shaped initial hypersurfaces $M_0$ in $\R$, if  $0<p<1$, and for strictly convex $M_0$ in case $1<p<\un$, where the curvature function $F$ is supposed to be homogeneous of degree 1, monotone and concave. The initial hypersurface is also assumed to be \tit{admissible}, i.e., its principal curvatures lie in the interior of the domain of $F$. 

Our results can be summarized in  two theorems:
\bt\lat{1.1}
Let $\C\su\R[n]$ be an open,  convex, symmetric cone containing the positive cone $\C_+$, and let $F\in C^{m,\al}(\C)\ii C^0(\bar\C)$, $4\le m\le \un$, $0<\al<1$, be a symmetric, strictly monotone and concave curvature function homogeneous of degree $1$ satisfying
\begin{equation}\lae{1.4}
\fv F\C>0\q\wed\q \fv F{\pa\C}=0,
\end{equation}
i.e., assume $\C$ to be a defining cone for $F$. Let $M_0\su\R$, $n\ge 2$, be a closed, admissible, star-shaped hypersurface of class $C^{m+2,\al}$, then the flow equation \re{1.2}, where $0<p<1$, with initial hypersurface $M_0$ has a solution $x=x(t,\xi)$ which exists for all time. The leaves $M(t)$ can be written as graphs of a function $u=u(t,\xi)$ over $\Ss[n]$ such that
\begin{equation}
u\in H^{m+2+\al,\frac{m+2+\al}2}(\bar Q),
\end{equation}
where $Q$ is the cylinder
\begin{equation}
Q=[0,\un)\times \Ss[n]
\end{equation}
and the parabolic H\"older space is defined in the usual way, see e.g., \cite[Note 2.5.4]{cg:cp}. The flow converges to infinity and the properly rescaled leaves converge in $C^{m+2}(\Ss[n])$ to the unit sphere.
\et
\bt\lat{1.2}
Let $1<p<\un$ and suppose that the curvature function $F$ satisfies the assumptions of the preceding theorem, and assume furthermore that $\C=\C_+$. Then the flow equation \re{1.2} with initial hypersurface $M_0\in C^{m+2,\al}$, where $M_0$ is closed and strictly convex, has a solution $x=x(t,\xi)$ which is defined on a maximal finite interval $[0,T^*)$. The leaves $M(t)$ can be written as graphs of a function $u=u(t,\xi)$ over $\Ss[n]$ such that
\begin{equation}
\lim_{t\ra T^*}\inf_{\Ss[n]}u(t,\cdot)=\un
\end{equation}
and the properly rescaled leaves converge in $C^{m+2}(\Ss[n])$ to the unit sphere.
\et

\section{Definitions and Conventions}

The main objective of this section is to state the equations of Gau{\ss}, Codazzi,
and Weingarten for hypersurfaces $M$ in a Riemannian \di{(n+1)} manifold $N$. 
Geometric quantities in $N$ will be denoted by
$(\bar g_{\alpha\beta})$, $(\riema \alpha\beta\gamma\delta)$, etc., and those in $M$ by $(g_{ij}), (\riem
ijkl)$, etc. Greek indices range from $0$ to $n$ and Latin from $1$ to $n$; the
summation convention is always used. Generic coordinate systems in $N$ resp.
$M$ will be denoted by $(x^\alpha)$ resp. $(\x^i)$. Covariant differentiation will
simply be indicated by indices, only in case of possible ambiguity they will be
preceded by a semicolon, i.e., for a function $u$ in $N$, $(u_\alpha)$ will be the
gradient and
$(u_{\alpha\beta})$ the Hessian, but e.g., the covariant derivative of the curvature
tensor will be abbreviated by $\riema \alpha\beta\gamma{\delta;\e}$. We also point out that
\begin{equation}
\riema \alpha\beta\gamma{\delta;i}=\riema \alpha\beta\gamma{\delta;\e}x_i^\e
\end{equation}
with obvious generalizations to other quantities.

In local coordinates, $(x^\alpha)$ and $(\x^i)$, the geometric quantities of the
 hypersurface $M$ are connected through the following equations
\begin{equation}\lae{2.3}
x_{ij}^\alpha=- h_{ij}\n^\alpha
\end{equation}
the so-called \tit{Gau{\ss} formula}. Here, and also in the sequel, a covariant
derivative is always a \tit{full} tensor, i.e.,

\begin{equation}
x_{ij}^\alpha=x_{,ij}^\alpha-\ch ijk x_k^\alpha+\cha \beta\gamma\alpha x_i^\beta x_j^\gamma.
\end{equation}
The comma indicates ordinary partial derivatives.

In this implicit definition the \tit{second fundamental form} $(h_{ij})$ is taken
with respect to $-\n$.

The second equation is the \tit{Weingarten equation}
\begin{equation}
\n_i^\alpha=h_i^k x_k^\alpha,
\end{equation}
where we remember that $\n_i^\alpha$ is a full tensor.

Finally, we have the \tit{Codazzi equation}
\begin{equation}
h_{ij;k}-h_{ik;j}=\riema\alpha\beta\gamma\delta\n^\alpha x_i^\beta x_j^\gamma x_k^\delta
\end{equation}
and the \tit{Gau{\ss} equation}
\begin{equation}
\riem ijkl= \{h_{ik}h_{jl}-h_{il}h_{jk}\} + \riema \alpha\beta\gamma\delta x_i^\alpha x_j^\beta x_k^\gamma
x_l^\delta.
\end{equation}

Now, let us assume that
$N$ is a topological product $\R[]\times \mc S_0$, where $\mc S_0$ is a
compact Riemannian manifold, and that there exists a Gaussian coordinate system
$(x^\alpha)$, such that the metric in $N$ has the form 
\begin{equation}\lae{2.8b}
d\bar s_N^2=e^{2\psi}\{{dx^0}^2+\sigma_{ij}(x^0,x)dx^idx^j\},
\end{equation}
where $\sigma_{ij}$ is a Riemannian metric, $\psi$ a function on $N$, and $x$ an
abbreviation for the  components $(x^i)$,

Let $M=\graph \fv u\so$ be a spacelike hypersurface
\begin{equation}
M=\set{(x^0,x)}{x^0=u(x),\,x\in\mc S_0},
\end{equation}
then the induced metric has the form
\begin{equation}
g_{ij}=e^{2\psi}\{ u_iu_j+\sigma_{ij}\}
\end{equation}
where $\sigma_{ij}$ is evaluated at $(u,x)$, and its inverse $(g^{ij})=(g_{ij})^{-1}$ can
be expressed as
\begin{equation}\lae{2.10}
g^{ij}=e^{-2\psi}\{\sigma^{ij}-\frac{u^i}{v}\frac{u^j}{v}\},
\end{equation}
where $(\sigma^{ij})=(\sigma_{ij})^{-1}$ and
\begin{equation}\lae{2.11}
\begin{aligned}
u^i&=\sigma^{ij}u_j\\
v^2&=1+ \sigma^{ij}u_iu_j\equiv 1+\abs{Du}^2.
\end{aligned}
\end{equation}

The contravariant form of a normal vector of a graph looks like
\begin{equation}
(\n^\alpha)=\pm v^{-1}e^{-\psi}(1, -u^i).
\end{equation}

In the Gau{\ss} formula \re{2.3} we are free to choose any of two normals, but we stipulate that in general we use \begin{equation}
(\n^\alpha)= v^{-1}e^{-\psi}(1, -u^i).
\end{equation}
as normal vector. 

Look at the component $\alpha=0$ in \re{2.3}, then we obtain 
\begin{equation}\lae{2.16}
e^{-\psi}v^{-1}h_{ij}=-u_{ij}-\cha 000\mspace{1mu}u_iu_j-\cha 0i0
\mspace{1mu}u_j-\cha 0j0\mspace{1mu}u_i-\cha ij0.
\end{equation}
Here, the covariant derivatives a taken with respect to the induced metric of
$M$, and
\begin{equation}\lae{2.18}
-\cha ij0=e^{-\psi}\bar h_{ij},
\end{equation}
where $(\bar h_{ij})$ is the second fundamental form of the hypersurfaces
$\{x^0=\const\}$.
 
\section{First estimates}

Let $F\in C^{m,\al}(\C)\ii C^0(\bar\C)$, $m\ge 4$, $0<\al<1$, be a monotone and concave curvature function homogeneous of degree 1 and normalized such that
\begin{equation}
F(1,\dots,1)=n.
\end{equation}

We first look at the flow of geodesic spheres. Fix a point $p_0\in\R$ and consider polar coordinates with center $p_0$. Then the Euclidean metric can be expressed as
\begin{equation}
d\bar s^2=dr^2+r^2\s_{ij}dx^idx^j,
\end{equation}
where $\s_{ij}$ is the standard metric of the sphere $\Ss[n]$.

Spheres with center $p_0$ and radius $r$ are umbilical, their second fundamental form is given by
\begin{equation}
\bar h_{ij}=r^{-1}\bar g_{ij}.
\end{equation}
Hence, the flow equation
\begin{equation}\lae{3.4}
\dot x=-\F\nu,
\end{equation}
where $\F(F)=-F^{-p}$, can be reduced to
\begin{equation}\lae{3.5}
\dot r=\frac1{(nr^{-1})^p}=n^{-p}r^p,
\end{equation}
and, thus, in case $p\not=1$,
\begin{equation}
r=\{\tfrac{1-p}{n^p}t+r_0^{1-p}\}^\frac1{1-p},
\end{equation}
where
\begin{equation}
r(0)=r_0,
\end{equation}
and we conclude:
\br
If the initial hypersurface is a sphere, the flow \re{3.5} exists for all time, if $0<p<1$, and converges to infinity, while in case $p>1$, the flow blows up in finite time
\begin{equation}
T^*=\frac{n^p}{p-1}r_0^{1-p}.
\end{equation}
\er
As a corollary we obtain:
\bc
Let $M_0=\graph u_0$ be star-shaped and let $x=x(t,\xi)$ be a solution of the flow \re{3.4} and define $u(t,\xi)$ by  $M(t)=\graph u(t)$. Let $r_1, r_2$ be positive constants such that
\begin{equation}
r_1<u_0(\xi)<r_2\qq\A\,\xi\in\Ss[n],
\end{equation}
then $u(t)$ satisfies the estimates
\begin{equation}\lae{3.10}
\Theta(t,r_1)<u(t,\xi)<\Theta(t,r_2)\qq\A\, 0\le t<\min\{T^*, T^*(r_1), T^*(r_2)\},
\end{equation}
where
\begin{equation}
\Theta(t,r)=\{\tfrac{1-p}{n^p}t+r^{1-p}\}^\frac1{1-p}
\end{equation}
and where  $T^*(r_i)$ indicates the maximal time for which  the spherical flow with initial sphere of radius $r_i$ will exist.
\ec
\bp
The spheres with radii $\Theta(r,r_i)$ are the spherical solutions of the flow \re{3.4}  with initial spheres of radius $r_i$.

The flow $x=x(t,\xi)$ also satisfies a scalar flow equation
\begin{equation}\lae{3.12} 
\dot u=\frac1{vF^p},
\end{equation}
where the dot indicates the total time derivative, or,
\begin{equation}\lae{3.13}
\pde ut=\frac v{F^p}
\end{equation}
when we consider the partial time derivative.

Hence, the result is due to the maximum principle, since these are parabolic equations.
\ep
\bl\lal{3.3}
Let $0<p<1$ and $r_1\le r\le r_2$, then there are positive constants $c_1$, $c_2$ depending only on $r_1$, $r_2$ and $p$ such that
\begin{equation}\lae{3.14}
0<c_1\le u(t)\Theta^{-1}(t,r)\le c_2\qq\A\,0\le t<T^*,
\end{equation}
and the flow is compactly contained in $\R[n]$ for finite $t$.
\el
\bp
Obvious.
\ep
\bl
Let $1<p$, then the flow \re{3.4}  only exists in a finite time interval $[0,T^*)$ and there holds 
\begin{equation}\lae{3.15}
\limsup_{t\ra T^*}\max_{\Ss[n]}u(t,\cdot)=\un.
\end{equation}
\el

\bp
Recall that for $p>1$ the hypersurfaces are convex. In view of the estimates \re{3.10} the maximal time $T^*$ has to be finite in this case.

As we shall prove later in \frt{4.1}  the flow will remain smooth with uniform estimates as long as it stays in a compact domain, hence \re{3.15} must be valid.  
\ep
Let $0<r_0$ be such that for the function $\Theta(t,r_0)$, where $p>1$, the singularity
\begin{equation}
T^*(r_0)=T^*,
\end{equation} 
then we can prove:
\bl\lal{3.5}
Let $u$ be the solution of the scalar flow equation \re{3.13} and assume $p>1$ and that \re{3.15} is valid. Then there exists a positive constant $c$ such that
\begin{equation}\lae{3.17}
u(t,\xi)-c\le \Theta(t,r_0)\le u(t,\xi)+c\qq\A\,\xi\in \Ss[n],
\end{equation}
hence
\begin{equation}\lae{3.18}
\lim_{t\ra T^*}u(t,\xi)\Theta^{-1}(t,r_0)=1\qq\A\,\xi\in \Ss[n].
\end{equation}
\el
\bp
Assume without loss of generality that the origin is inside the convex body defined by $M_0$. Then the support function
\begin{equation}
\bar u=\spd x\nu
\end{equation}
of the flow hypersurfaces can be looked at as being defined on the Gau{\ss} image of $M(t)$ and $\bar u$ satisfies the parabolic equation
\begin{equation}\lae{3.20} 
\begin{aligned}
\dot{\bar u}=\pde{\bar u}t&=\tilde F^p(\bar u_{ij}+\bar u\s_{ij})\\
&=\tilde F(h_{ij},\s_{ij})=\tilde F^p(\ka_i^{-1})
\end{aligned}
\end{equation}
on $\Ss[n]$, where $\ka_i$ are the principal curvatures of $M(t)$, $h_{ij}$ the second fundamental form and $\tilde F$ the inverse curvature function of $F$ which is defined by
\begin{equation}
\tilde F(\ka_i)=\frac1{F(\ka_i^{-1})}\qq\A\,(\ka_i)\in \C_+.
\end{equation}
Equation \re{3.20} can be easily derived, see e.g., \cite[Section 2]{urbas:convex}.

In view of the results in \cite[Theorem 3.1]{chow-gulliver} and \cite[Theorem 3.1]{mccoy:vol1} the solution $\bar u$ of \re{3.20} satisfies the a priori estimate
\begin{equation}\lae{3.21} 
\osc\bar u\le c,
\end{equation}
where $c=c(\bar u_0)$, from which we immediately deduce
\begin{equation}
\osc u\le c,
\end{equation}
since $u(t,\xi)=\bar u(t,\xi)$, when $\xi$ is an extremal point. The estimate \re{3.17} then is due to the fact that for any $t\in [0,T^*)$ there exists $\xi_t$ such that
\begin{equation}
u(t,\xi_t)=\Theta(t,r_0),
\end{equation}
\cf the arguments in the proof of \cite[Lemma 5.1]{oliver:2004surfaces}. 
\ep

Next, we want to prove a priori estimates for $v$, or equivalently, for
\begin{equation}
\abs{Du}^2=u^{-2}\s^{ij}u_iu_j=\s^{ij}\f_i\f_j\equiv\abs{D\f}^2,
\end{equation}
where
\begin{equation}
\f=\log u.
\end{equation}

Let us consider instead of the Euclidean metric  a more general metric
\begin{equation}
d\bar s^2=dr^2+\vt(r)^2\s_{ij}dx^idx^j.
\end{equation}
The second fundamental form can then be expressed as
\begin{equation}\lae{3.26}
h_{ij}v^{-1}=-u_{ij}+\bar h_{ij}=-u_{ij}+\dot\vt\vt\s_{ij},
\end{equation}
where the covariant derivatives are taken with respect to the induced metric.

Define the metric
\begin{equation}
\tilde \s_{ij}=\vt^2(u)\s_{ij},
\end{equation}
and denote covariant differentiation with respect to this metric by a semi-colon, then
\begin{equation}
h_{ij}v^{-1}=-v^{-2}u_{;ij}+\dot\vt\vt \s_{ij},
\end{equation}
\cf \cite[Lemma 2.7.6]{cg:cp}, and we conclude further
\begin{equation}
\begin{aligned}
h^i_j&=g^{ik}h_{kj}\\
&=v^{-1}\vt^{-1}\{-(\s^{ik}-v^{-2}\f^i\f^k)\f_{jk}+\dot\vt \de^i_j\},
\end{aligned}
\end{equation}
where $\s^{ij}$ is the inverse of $\s_{ij}$,
\begin{equation}
\f=\int_{r_0}^u\vt^{-1},
\end{equation}
\begin{equation}
\f^i=\s^{ik}\f_k,
\end{equation}
and $\f_{jk}$ are the second covariant derivatives of $\f$ with respect to the metric $\s_{ij}$.

Thus, the scalar curvature equation \re{3.13} can now be expressed as
\begin{equation}
\dot u=\frac v{F^p},
\end{equation}
or equivalently,
\begin{equation}\lae{3.30} 
\dot\f=\vt^{-1}\dot u =\frac {\vt^{p-1}v}{F^p(\vt h^i_j)}\equiv \frac{\vt^{p-1}v}{F^p(\tilde h^i_j)},
\end{equation}
where
\begin{equation}
\tilde h^i_j=v^{-1}\{-(\s^{ik}-v^{-2}\f^i\f^k)\f_{jk}+\dot\vt \de^i_j\}.
\end{equation}
Let
\begin{equation}\lae{3.35} 
\tilde g_{ij}=\f_i\f_j+\s_{ij},
\end{equation}
then we consider the eigenvalues of
\begin{equation}
\tilde h_{ij}=\tilde g_{ik}\tilde h^k_j
\end{equation}
with respect to this metric and we define $F^{ij}$ \resp $F^i_j$ accordingly
\begin{equation}\lae{3.37}
F^{ij}=\pdc F{\tilde h}{{ij}}
\end{equation}
and
\begin{equation}
F^i_j=\pdm F{\tilde h}ij=\tilde g_{jk}F^{ik}.
\end{equation}
Note that $\tilde h_{ij}$ is symmetric, since $h_{ij}$ and $\tilde g_{ij}$ can be diagonalized simultaneously. We also emphasize that
\begin{equation}\lae{3.36}
\abs{Du}^2=\s^{ij}\f_i\f_j\equiv\abs{D\f}^2.
\end{equation}
\bl\lal{3.6}
Let $u$ be a solution of the scalar curvature equation
\begin{equation}
\dot u=\frac v{F^p},
\end{equation}
where $0<p<1$, $\dot\vt\ge 0$ and $\Ddot\vt\ge0$, then
\begin{equation}\lae{3.41}
\abs{Du}^2\le \sup_{\Ss[n]}\abs{Du_0}^2
\end{equation}
during the evolution.

Moreover, if $\vt(u)=u$ and
\begin{equation}\lae{3.42}
F(\tilde h^i_j)=uF(h^i_j)\ge c_0>0,
\end{equation}
then
\begin{equation}\lae{3.43}
\abs{Du}^2\le \sup_{\Ss[n]}\abs{Du_0}^2\frac{b^\ga}{(at+b)^\ga}
\end{equation}
with positive constants $a$, $b$, and $\ga$.
\el
\bp
\cq{\re{3.41}}\q In view of  \re{3.36}, we may estimate 
\begin{equation}
w=\tfrac12\abs{D\f}^2.
\end{equation}
Differentiating equation \re{3.30} covariantly with respect to
\begin{equation}
\f^kD_k
\end{equation}
we deduce 
\begin{equation}\lae{3.46}
\begin{aligned}
\dot w&=(p-1)\vt^{p-1}\dot\vt F^{-p}\abs{D\f}^2v+\vt^{p-1}v_k\f^kF^{-p}\\
&\q +p\vt^{p-1}vF^{-(p+1)}\{v^{-1}F^j_i\tilde h^i_jv_k\f^k+ v^{-1} F^k_l\tilde g^{lr}w_{kr} \\
&\q-v^{-1}F^k_l\tilde g^{lr}\f_{ik}\f^i_r
+v^{-1}F^k_l\tilde g^{lr}_{\hp{lr};i}\f^i\f_{kr}+v^{-1}F^k_l\tilde g^{lr}\f_r\f_k\\
&\q-v^{-1}F^k_l\tilde g^{lr}\s_{kr}\abs{D\f}^2
 -2v^{-1}F^k_k\Ddot\vt\vt w\},
\end{aligned}
\end{equation}
where covariant derivatives with respect to the metric $\s_{ij}$ are simply denoted by indices, if no ambiguities are possible, and by a semi-colon otherwise. In deriving the previous equation we also used the Ricci identities and the properties of the Riemann curvature tensor of $\Ss[n]$.

Now, let $0<T<T^*$ and suppose that 
\begin{equation}
\sup_{\Q_T}w,\qq Q_T=[0,T]\times \Ss[n],
\end{equation}
is attained at $(t_0,x_0)$ with $t_0>0$. Then the maximum principle implies
\begin{equation}\lae{3.45}
\begin{aligned}
0&\le -F^k_l\tilde g^{lr}\f_{ik}\f^i_r+(F^k_l\tilde g^{lr}\f_r\f_k-F^k_l\tilde g^{lr}\s_{kr}\abs{D\f}^2)\\
&\q -2F^k_k\Ddot\vt\vt w.
\end{aligned}
\end{equation}
The right-hand side, however, is strictly negative, if $w>0$, hence $t_0>0$ is not possible, since we didn't assume $M_0$ to be a sphere, and we conclude
\begin{equation}
w\le \sup_{\Ss[n]}w(0).
\end{equation} 

\cq{\re{3.43}}\q Define
\begin{equation}
\bar w=\sup_{\Ss[n]}w(t,\cdot)=w(t,\xi_t),
\end{equation}
then $\bar w$ is Lipschitz and for almost every $t\in[0,T^*)$ there holds
\begin{equation}
\dot{\bar w}=\pde {w(t,\xi_t)}t.
\end{equation}
Applying the maximum principle and the definition of $\vt$ we infer from \re{3.46}
\begin{equation}
\dot{\bar w}\le2(p-1)u^{p-1}F^{-p}\bar w v,
\end{equation}
where $F=F(\tilde h^i_j)$. Using  the notations and results of \rl{3.3} we further deduce
\begin{equation}
\dot{\bar w}\le (p-1)2F^{-p}\tilde u^{p-1}\bar w\Theta^{p-1}v,
\end{equation}
where $\tilde u=u\Theta^{-1}$, hence
\begin{equation}
(\log\bar w)'\le -c\frac1{at+b}
\end{equation}
in view of the assumption \re{3.42} and the estimate \re{3.41}; $a$,  $b$ and $c$ are positive constants. The relation \re{3.43} is then an immediate consequence with $\ga=\frac ca$.
\ep
A similar a priori estimate is also valid in case $p>1$.
\bl
Let $p>1$ and assume \re{3.15} to be satisfied, then
\begin{equation}
v-1\le c\,\Theta^{-1},
\end{equation}
where $\Theta=\Theta(t,r_0)$ as in \re{3.17}, i.e.,
\begin{equation}
\lim_{t\ra T^*}\norm{Du}=0.
\end{equation}
\el
 \bp
 Let $\bar u$ be the support function and let
 \begin{equation}
\tilde  u=u\Theta^{-1}\q\wed\q\tilde{\bar u}=\bar u\Theta^{-1},
\end{equation}
then we deduce from \re{3.17} and \re{3.21} and the relation
\begin{equation}
\bar u=\spd x\nu,
\end{equation}
or equivalently,
\begin{equation}
v=u\bar u^{-1},
\end{equation}
that
\begin{equation}
\begin{aligned}
v-1&=(u-\bar u)\bar u^{-1}=(\tilde u-\tilde{\bar u})\tilde{\bar u}^{-1}\\
&=(\tilde u-1)\tilde{\bar u}^{-1}+(1-\tilde{\bar u})\tilde{\bar u}^{-1}\\
&\le c\,\Theta^{-1}.
\end{aligned}
\end{equation}
 \ep

Next let us state the evolution equations for the geometric quantities $\F$, $u$, and $\chi$, where
\begin{equation}
\chi=\spd x\nu^{-1}=vu^{-1}.
\end{equation}
The term $\F=-F^{-p}$ satisfies the evolution equation
\begin{equation}\lae{3.62}
\F'-\dot\F F^{ij}\F_{ij}=\dot\F F^{ij}h_{ik}h^k_j\F,
\end{equation}
\cf \cite[Lemma 2.4.8]{cg:cp}, where we use the notation
\begin{equation}
\F'=\df {(\F\circ F)}t\q\wed\q \dot\F=\df\F r,
\end{equation}
here $\F=\F(r)$; for $u$ we have 
\begin{equation}\lae{3.63}
\dot u-\dot\F F^{ij}u_{ij}=v^{-1}(p+1)F^{-p}-\dot\F F^{ij}\bar h_{ij},
\end{equation}
which can be easily deduced from the scalar curvature equation \re{3.12} and the expression \re{3.26} for the second fundamental form in view of the homogeneity of $F$. Finally, $\chi$ satisfies
\begin{equation}\lae{3.64}
\dot\chi-\dot\F F^{ij}\chi_{ij}=-\dot\F F^{ij}h_{ik}h^k_j\chi-2\chi^{-1}\dot\F F^{ij}\chi_i\chi_j+\{\dot\F F+\F\}\frac{\bar H}nv\chi,
\end{equation}
\cf \cite[Lemma 5.8]{cg:spaceform}.

The spherical solution $\Theta$, which we also call the scale factor, satisfies
\begin{equation}
\dot\Theta=n^{-p}\Theta^p
\end{equation}
for $p\not=1$, hence, the scaled quantities
\begin{equation}
\tilde\F=\F\Theta^{-p}\q\wed\q\tilde u=u\Theta^{-1}\q\wed\q\tilde\chi=\chi\Theta
\end{equation}
satisfy similar equations with an additional summand at the right-hand side. Let us only state the equation for $\tilde u$ explicitly
\begin{equation}\lae{3.67}
\begin{aligned}
\dot{\tilde u}-\dot\F F^{ij}\tilde u_{ij}&=v^{-1}(p+1)\tilde F^{-p}\Theta^{p-1}-p\tilde F^{-(p+1)}F^{ij}\bar g_{ij}\tilde u^{-1}\Theta^{p-1}\\
&\qq-n^{-p}\tilde u\Theta^{p-1}.
\end{aligned}
\end{equation}
\br\lar{3.8}
The previous equation---and also any other equation satisfied by  a properly rescaled geometric quantity of the flow \re{3.4}---has a homogeneity of order $(p-1)$ in $\Theta$ of the right-hand side. Of course the elliptic part of the differential operator shows the same homogeneity though it is hidden at the moment. Indeed the elliptic part of equation \re{3.67} can be expressed in the form
\begin{equation}\lae{3.68}
\begin{aligned}
-\dot\F F^{ij}\tilde u_{ij}&=-\dot{\tilde\F}\tilde u^{-2}\{F^{ij}\tilde u_{ij}-2F^{ij}\tilde u_i\f_j+\tilde u_k\f^kF^{ij}\tilde g_{ij}\}\Theta^{p-1}
\end{aligned}
\end{equation}
where the spatial covariant derivatives on the right-hand side are now taken with respect to the metric $\tilde g_{ij}$ in \re{3.35}, $\f=\log u$, and $F^{ij}$ is defined as in \re{3.37}. This observation combined with a simple variable transformation with respect to time will later enable us to immediately deduce higher order estimates once we have proved $C^2$-estimates and compactness of the rescaled principal curvatures in the defining cone $\C$, in case $0<p<1$, \resp in $\C_+$ for $p>1$.

Note that 
\begin{equation}
F^{ij} g_{ij}=F^i_i
\end{equation}
is scaling invariant and will satisfy
\begin{equation}
n\le F^i_i\le c,
\end{equation}
when the rescaled principal curvatures
\begin{equation}
\ka_i\Theta
\end{equation}
are compactly contained in the respective cones. 
\er

We can now prove that the assumption \re{3.42} is satisfied.
\bl\lal{3.9}
Let $0<p<1$, then the assumption \re{3.42} is valid, i.e., there exists a constant $c_0$ such that
\begin{equation}\lae{3.72}
F(\tilde h^i_j)\ge c_0>0,
\end{equation}
where
\begin{equation}
\tilde h^i_j=h^i_j\Theta=h^i_ju\tilde u^{-1}.
\end{equation}
\el
\bp
Define
\begin{equation}
w=\log(-\tilde\F)+\log\tilde\chi
\end{equation}
and let $0<T<T^*$ be arbitrary. Applying the maximum principle to $w$ in the cylinder $Q_T$ and assuming that the maximum is attained at a time $t_0>0$ we infer from \re{3.62} and \re{3.64} that
\begin{equation}
0\le (p-1)\tilde F^{-p}\tilde u^{-1}v\Theta^{p-1}+n^{-p}(1-p)\Theta^{p-1}
\end{equation}
implying \re{3.72}. 
\ep

The same result is also valid in case $p>1$.
\bl\lal{3.10}
Let $p>1$ and assume the relation \re{3.15}, then there exists a constant $c_0$ such that
\begin{equation}
F(\tilde h^i_j)\ge c_0>0.
\end{equation}
\el
\bp
Let $0<T<T^*$ be arbitrary. Looking at a supremum in $Q_T$ for the function
\begin{equation}
w=\log(-\tilde\F)+\log\tilde\chi+\tilde u
\end{equation}
and applying the maximum principle we obtain
\begin{equation}
\begin{aligned}
0&\le \dot\F F^{ij}(\log(-\tilde\F))_i(\log(-\tilde\F))_j-\dot\F F^{ij}(\log\tilde\chi)_i(\log\tilde\chi)_j\\
&\q+c(p-1)\tilde F^{-p}\Theta^{p-1}-p\tilde F^{-(p+1)}F^{ij}\bar g_{ij}\tilde u^{-1}\Theta^{p-1}.
\end{aligned}
\end{equation}
Ignoring the derivatives for the moment we see that the last term is dominating with the right sign.

To estimate the derivatives we use the fact that $w_i=0$ and deduce
\begin{equation}\lae{3.79}
\begin{aligned}
&\dot\F F^{ij}(\log(-\tilde\F))_i(\log(-\tilde\F))_j-\dot\F F^{ij}(\log\tilde\chi)_i(\log\tilde\chi)_j=\\
&\qq\dot\F F^{ij}\tilde u_i\tilde u_j +2\dot\F F^{ij}(\log\tilde\chi)_i\tilde u_j.
\end{aligned}
\end{equation}
The first term on the right-hand side is of the order $\tilde F^{-(p+1)}$, but $\norm{D\tilde u}$ vanishes if $t$ tends to $T^*$, and we shall show in \frl{4.4} that $\tilde F^{-1}$ is bounded when $t$ stays in compact subsets of $[0,T^*)$.

The mixed term on the right-hand side of \re{3.79} is nonpositive since the leaves $M(t)$ are supposed to be strictly convex and $\chi$ is equal to
\begin{equation}
\chi=\spd x\nu^{-1};
\end{equation}
the claim then follows by using the Weingarten equation. Hence, $w$ is a priori bounded from above and the lemma proved.
\ep

\section{$C^2$-estimates and maximal existence}
In this section we shall prove uniform estimates for the rescaled second fundamental
\begin{equation}
\tilde h^i_j=h^i_j\Theta^{-1},
\end{equation}
that in case $0<p<1$ the flow exists for all time, and that in case $p>1$ the maximal time $T^*$ is indeed characterized by a blow up of the flow.

Let us start with the latter result.
\bt\lat{4.1}
Let $p>1$ and
let the initial hypersurface $M_0\in C^{m+2,\al}$, $4\le m\le\un$, $0<\al<1$, be strictly convex,  then the solution of the curvature flow
\begin{equation}
\dot x=-\F\nu
\end{equation}
exists for $[0,T^*)$ and belongs to the parabolic H\"older space $H^{m+2+\al,\frac{m+2+\al}2}(Q)$, where
\begin{equation}
Q=[0,T^*)\times \Ss[n].
\end{equation}
The flow \fre{3.4} satisfies uniform estimates in this function class as long as it stays in a compact subset  $\bar\Om\su \R$. Moreover, the principal curvatures are strictly convex
\begin{equation}\lae{4.2}
0<c_1\le \ka_i\le c_2\qq\A\,1\le i\le n,
\end{equation}
where the constants $c_1$, $c_2$ depend on $\bar\Om$, $p$, and $M_0$. Hence, the singularity $T^*$ is characterized to be the blow up time of the flow, i.e., the relation \fre{3.15} is valid.
\et

For the proof we need several lemmata.
\bl\lal{4.2}
Let $M$ be a closed, convex hypersurface which is represented as the graph of a $C^1$-function $u$ over $\Ss[n]$. Assume that $u$ is bounded by
\begin{equation}
0<r_1\le u\le r_2,
\end{equation}
then
\begin{equation}
v\le c(r_1,r_2).
\end{equation}
\el
The lemma is proved in \cite[Theorem 2.7.10]{cg:cp}.
\bl\lal{4.3}
Let $M(t)$ be a solution of the flow equation \fre{3.4} for any given $0<p$, then
\begin{equation}\lae{4.5}
F\le \sup_{M_0}F
\end{equation}
\el
\bp
Use equation \fre{3.62} and apply the parabolic maximum principle.
\ep
\bl\lal{4.4}
Let $p>1$ and let $M(t)$ be a solution of the flow \re{3.4}. Let $0<T<T^*$ be arbitrary and assume
\begin{equation}\lae{4.6}
u\le r_2\qq\A\, 0\le t\le T,
\end{equation}
then there exists a constant $c_0$ such that
\begin{equation}\lae{4.7}
0<c_0\le F\qq\A\,0\le t\le T,
\end{equation}
where $c_0=c_0(r_2,p,M_0)$ is independent of $T$.
\el
\bp
The proof is almost identical to the proof of \frl{3.10}, the only difference is that this time we omit the scale factor and apply the maximum principle to the function
\begin{equation}
w=\log(-\F)+\log\chi+u.
\end{equation}
\ep
\bl\lal{4.5} 
Let $p>1$ and $0<T<T^*$ be arbitrary; assume that the flow satisfies the estimate \re{4.6}, then there exists a constant $c$ depending on $r_2$ and $M_0$, but not on $T$, such that the principal curvatures $\ka_i$ satisfy
\begin{equation}
\ka_i\le c.
\end{equation}
\el
\bp
The mixed tensor $h^i_j$ satisfies the evolution equation 
\begin{equation}\lae{4.10}
\begin{aligned}
\dot h^i_j-\dot\F F^{kl}h^i_{j;kl}&=\dot\F F^{kl}h_{kr}h^r_kh^i_j+(\F-\dot\F F)h^{ki}h_{kj}\\
&\q +\Ddot\F F_jF^i
+\dot\F F^{kl,rs}h_{kl;j}h_{rs;}^{\hp{rs;}i} 
\end{aligned}
\end{equation}
\cf \cite[Lemma 2.4.3]{cg:cp}.

Define $\zeta$ and $w$ by
\begin{equation}
\zeta=\sup\set{h_{ij}\h^i\h^j}{\norm{\h}=1}
\end{equation}
and
\begin{equation}\lae{4.12}
w=\log\zeta+\log\chi.
\end{equation}
We claim that $w$ is bounded from above in $Q_T$ by a constant independent of $T$.

Let $x_0=x_0(t_0,\xi_0)$ with $0<t_0\le T$ be a point in $M(t_0)$ such that
\begin{equation}
\sup_{M_0}w<\sup\set {\sup_{M(t)} w}{0<t\le T}=w(x_0).
\end{equation}

We then introduce a Riemannian normal coordinate system $(\x^i)$ at $x_0\in
M(t_0)$ such that at $x_0=x(t_0,\x_0)$ we have
\begin{equation}
g_{ij}=\delta_{ij}\q \tup{and}\q \zeta=h_n^n.
\end{equation}

Let $\tilde \h=(\tilde \h^i)$ be the contravariant vector field defined by 
\begin{equation}
\tilde \h=(0,\dotsc,0,1),
\end{equation}
and set
\begin{equation}
\tilde \zeta=\frac{h_{ij}\tilde \h^i\tilde \h^j}{g_{ij}\tilde \h^i\tilde \h^j}\raise 2pt
\hbox{.}
\end{equation}

$\tilde\zeta$ is well defined in neighbourhood of $(t_0,\x_0)$.

Now, define $\tilde w$ by replacing $\zeta$ by $\tilde \zeta$ in \re{4.12}; then, $\tilde w$
assumes its maximum at $(t_0,\x_0)$. Moreover, at $(t_0,\x_0)$ we have 
\begin{equation}
\dot{\tilde \zeta}=\dot h_n^n,
\end{equation}
and the spatial derivatives do also coincide; in short, at $(t_0,\x_0)$ $\tilde \zeta$
satisfies the same differential equation \re{4.10} as $h_n^n$. For the sake of
greater clarity, let us therefore treat $h_n^n$ like a scalar and pretend that $w$
is defined by 
\begin{equation}
w=\log h_n^n+ \log\chi.
\end{equation} 

From equations \re{4.10} and \fre{3.64} we then infer
\begin{equation}
0\le -(p+1)F^{-p}h^n_n+(p-1)cF^{-p},
\end{equation}
where we used the concavity of  $F$. The constant $c$ is equal to
\begin{equation}
\sup_{Q_T}u^{-1}v;
\end{equation}
hence, $h^n_n$ is a priori bounded independent of $T$.
\ep
It remains to prove a lower positive bound for the $\ka_i$, but this immediately follows from \re{4.7} and the assumptions \fre{1.4}. Hence, we can state:
\bl\lal{4.7}
Let $p>1$ and let $0<T<T^*$ be arbitrary; assume that the flow satisfies the estimate \re{4.6}, then there exists a positive constant $\ka_0$ depending on $r_2$, $T^*$, and $M_0$, but not on $T$, such that the principal curvatures are bounded from below by 
\begin{equation}
0<\ka_0\le\ka_i\qq\A\, 1\le i\le n.
\end{equation}
\el

We can now prove \rt{4.1}.
\bp[Proof of \rt{4.1}]\lapp{4.1}
In the preceding lemmata we proved uniform $C^2$-estimates as well as the  estimate \re{4.2} for convex solutions of the flow equations provided the flow stays in a compact subset of $\R$.

The corresponding scalar curvature equation \fre{3.13} is then a nonlinear parabolic equation where the elliptic part is concave, and, because of the a priori estimates, uniformly elliptic provided the flow is compactly contained in $\R$. Hence, we can apply the Krylov-Safonov estimates yielding uniform H\"older estimates for $\dot u$ and $D^2u$ estimates. Now, the linear theory and the parabolic Schauder estimates can be applied; for details see  e.g., \cite[Chapter 2.6]{cg:cp} and \cite[Section 6]{cg:survey}.
\ep

Let us now derive corresponding estimates for the rescaled quantities, where at the moment we still assume $p>1$.

We already know that
\begin{equation}
\tilde u=u\Theta^{-1}
\end{equation}
is uniformly bounded from above and against zero, \cf \frl{3.5}, that $\norm{D\tilde u}$ vanishes if $t$ tends to $T^*$ and that
\begin{equation}
\tilde\chi=\chi \Theta=vu^{-1}\Theta=v\tilde u^{-1}
\end{equation}
is uniformly bounded from below and from above
\begin{equation}
0<c_1\le\tilde\chi\le c_2
\end{equation}
such that
\begin{equation}
\lim_{t\ra T^*}\tilde\chi=1.
\end{equation}

Furthermore, we already established $F(\tilde h^i_j)\ge c_0$, \cf \frl{3.10}; hence, it remains to prove
\begin{equation}
0<c_1\le \tilde\ka_i\le c_2\qq\A\,1\le i\le n,
\end{equation}
where
\begin{equation}
\tilde\ka_i=\ka_i\Theta.
\end{equation}
Let us also recall the notation
\begin{equation}
\tilde F=F\Theta=F(\tilde h^i_j).
\end{equation}

\bl\lal{4.8}
Let $p>1$ then the rescaled principal curvatures $\tilde\ka_i$ of the flow hypersurfaces are uniformly bounded from above
\begin{equation}\lae{4.37}
\tilde\ka_i\le c\qq\A\, 1\le i\le n,
\end{equation}
during the evolution.
\el
\bp
The proof is similar to the proof of \rl{4.5}. We define
\begin{equation}\lae{4.38} 
\zeta=\sup\set{\tilde h_{ij}\h^i\h^j}{\norm\h=1}
\end{equation}
and
\begin{equation}
w=\log\zeta+\log\tilde\chi+\lam\tilde u,
\end{equation}
where $\lam>1$ is large.

Fix an arbitrary $0<T<T^*$ and assume that
\begin{equation}
\sup_{M_0}w<\sup\set{\sup_{M(t)}w}{0\le t\le T}=w(x_0)),
\end{equation}
where $x_0$=$x(t_0,\xi_0)$ such that $t_0>0$.

Arguing as in the proof of \rl{4.5} we may replace $\zeta$ by $\tilde h^n_n$; applying the maximum principle and observing \frr{3.8} then yields
\begin{equation}
\begin{aligned}
0&\le -(p+1)\tilde F^{-p}\tilde h^n_n\Theta^{p-1}+\lam c\tilde F^{-p}\Theta^{p-1}-(\lam\tilde u-2)n^{-p}\Theta^{p-1}\\
&\q -\lam c \tilde F^{-(p+1)}F^{ij}g_{ij}\Theta^{p-1}+p\tilde F^{-(p+1)}F^{ij}(\log\tilde h^n_n)_i(\log\tilde h^n_n)_j\Theta^{p-1}\\
&\q -p\tilde F^{-(p+1)}F^{ij}(\log\tilde\chi)_i(\log\tilde\chi)_j\Theta^{p-1}.
\end{aligned}
\end{equation}
To estimate the two terms involving the derivatives we use the fact that $Dw=0$ and deduce, abbreviating the difference of the two terms by $D$,
\begin{equation}\lae{4.42}
\begin{aligned}
D=p\lam^2\tilde F^{-(p+1)}F^{ij}\tilde u_i\tilde u_j+2p\tilde F^{-(p+1)}F^{ij}(\log\tilde\chi)_i\tilde u_j\Theta^{p-1}.
\end{aligned}
\end{equation}
The last term on the right-hand side is nonpositive, \cf the arguments after equation \fre{3.79}, while
\begin{equation}
F^{ij}\tilde u_i\tilde u_j\le F^{ij}g_{ij}\norm{D\tilde u}^2.
\end{equation}
Choosing now $T$ close to $T^*$ and assuming  $w(x_0)$ to be  large enough such that $t_0$ is close to $T$ and hence $\norm {D\tilde u}^2$ is very small, the first term on the right-hand side of  \re{4.42} can be absorbed by
\begin{equation}
-\lam c\tilde F^{-(p+1)}F^{ij}g_{ij}\Theta^{p-1}
\end{equation}
and we conclude that $w$ is a priori bounded.
\ep
\br
The estimate \re{4.37} implies that
\begin{equation}
\tilde F\le n c.
\end{equation}
\er
\bp
Follows immediately from the monotonicity and homogeneity of $F$.
\ep
Combining the results of \rl{4.8} and \frl{3.10} we deduce:
\bl
Let $p>1$, then the rescaled principal curvatures of the flow hypersurfaces are uniformly positive
\begin{equation}
0<c_1\le \tilde\ka_i\qq\A\,1\le i\le n.
\end{equation}
\el

Let us now prove corresponding estimates in case $0<p<1$.
\bl
Let $0<p<1$, then the principal curvatures of the flow hypersurfaces are uniformly bounded during the evolution
\begin{equation}
\ka_i\le c\qq\A\,1\le i\le n.
\end{equation}
\el
\bp
The proof is identical to the proof of \rl{4.5}; we then deduce
\begin{equation}
\ka_i\le c \sup_{M_0}\max_{1\le i\le n}\ka_i,
\end{equation}
since the maximum principle yields
\begin{equation}
0\le -(p+1)F^{-p}h^n_n+(p-1)F^{-p}u^{-1}v<0,
\end{equation}
if the maximum of $w$ is attained at a time $t_0>0$.
\ep
Combining this estimate with the already known estimates for $\abs{Du}$, $\tilde u$, and $F$, \cf \frl{3.3}, \frl{3.6}, and \frl{3.9}, we conclude that the flow exists for all time and converges to infinity, see the final argument in the proof of  \rt{4.1} on page \rpp{4.1}.

Next, let us estimate the rescaled principal curvatures.
\bl
Let $0<p<1$, then the rescaled principal curvatures of the flow hypersurfaces 
\begin{equation}
\tilde\ka_i=\ka_i\Theta
\end{equation}
are uniformly bounded from above
\begin{equation}
\tilde\ka_i\le c\qq\A\, 1\le i\le n.
\end{equation}
\el
\bp
The proof is similar to the proof of \rl{4.8}. Define $\zeta$ as in \re{4.38}, choose $0<T<\un$ very large, and suppose that
\begin{equation}
w=\log\zeta +\lam\log\tilde\chi+\mu\tilde u,
\end{equation}
where $\lam$ and $\mu=\mu(\lam)$ are large.

Applying the maximum principle in the cylinder $Q_T$ and replacing $\zeta$ by $\tilde h^n_n$ we obtain
\begin{equation}
\begin{aligned}
0&\le -(p+1)\tilde F^{-p}\tilde h^n_n\Theta^{p-1}-(\lam-1)p\tilde F^{-(p+1)}F^{ij}\tilde h_{ik}\tilde h^k_j\Theta^{p-1}\\
&\q+\lam c \tilde F^{-p}\Theta^{p-1}-(\mu\tilde u-\lam-1)n^{-p}\Theta^{p-1}\\
&\q-\mu v^{-2}p\tilde F^{-(p+1)}F^{ij}g_{ij}\tilde u^{-1}\Theta^{p-1} \\
&\q+ p\tilde F^{-(p+1)}F^{ij}(\log\tilde h^n_n)_i(\log\tilde h^n_n)_j\Theta^{p-1}\\
&\q-\lam p\tilde F^{-(p+1)}F^{ij}(\log\tilde\chi)_i(\log\tilde\chi)_j\Theta^{p-1}.
\end{aligned}
\end{equation}
The critical terms are those containing the derivatives. Let us abbreviate by $D$ the sum of the last two terms of the right-hand side of the preceding inequality. Using the fact that $Dw=0$ we conclude
\begin{equation}
\begin{aligned}
D&=(\lam^2-\lam)p\tilde F^{-(p+1)}F^{ij}(\log\tilde\chi)_i(\log\tilde\chi)_j\Theta^{p-1}\\
&\q+\mu^2p\tilde F^{-(p+1)}F^{ij}\tilde u_i\tilde u_j -2\lam\mu p\tilde F^{-(p+1)}F^{ij}(\log\tilde\chi)_i\tilde u_j\Theta^{p-1}.
\end{aligned}
\end{equation}
Now, there holds
\begin{equation}
\tilde\chi=\spd x\nu^{-1}\Theta
\end{equation}
and hence
\begin{equation}
\begin{aligned}
\tilde\chi_i&=-\chi^2 h^k_i\spd{x_k}x\Theta\\
&=-\chi^2 h^k_iu_ku\Theta\\
&=-\tilde\chi^2\tilde h^k_i\tilde u_k\tilde u,
\end{aligned}
\end{equation}
which implies
\begin{equation}
F^{ij}(\log\tilde\chi)_i(\log\tilde\chi)_j\le v^2F^{ij}\tilde h^k_i\tilde h_{kj}\norm{D\tilde u}^2.
\end{equation}
Thus, $D$ can be absorbed by
\begin{equation}
-(\lam-1)p\tilde F^{-(p+1)}F^{ij}\tilde h_{ik}\tilde h^k_j\Theta^{p-1}-\mu v^{-2}p\tilde F^{-(p+1)}F^{ij}g_{ij}\tilde u^{-1}\Theta^{p-1},
\end{equation}
if $t_0$ is large enough; hence $\tilde h^n_n$ is a priori bounded.
\ep
Combing this result with the estimate \fre{3.72} we deduce:
\bc
Let $0<p<1$, then the rescaled principal curvatures of the flow hypersurfaces stay in a compact subset of the defining cone $\C$. 
\ec

\section{Convergence of the rescaled flow}
In this section we do not need to distinguish between the cases $p<1$ and $p>1$. The results and the arguments are valid for any positive $p\not=1$.

In the previous sections we have proved uniform estimates for the rescaled flow of class $C^2$ and we also proved that the rescaled elliptic operator $F\Theta$ is uniformly elliptic, or more precisely:
\bl
The rescaled curvature function
\begin{equation}
F(\tilde h^i_j) =F(h^i_j)\Theta
\end{equation}
represents a uniformly elliptic differential operator for the function $\log\tilde u$, where
\begin{equation}
\tilde h^i_j=v\tilde u^{-1}\{-\f^i_j+v^{-2}\de^i_j\}, 
\end{equation}
\begin{equation}
\f=\log u,
\end{equation}
$\tilde g_{ij}$ is the metric defined in \fre{3.35}, and where the covariant derivatives of $\f$ are with respect to this metric. The rescaled scalar curvature equation \fre{3.13} then takes the form
\begin{equation}
\pde{\tilde u}t=vF^{-p}\Theta^{p-1}-n^{-p}\tilde u\Theta^{p-1},
\end{equation}
where $F=F(\tilde h^i_j)$.

Defining  $\tau=\tau(t)$ by the relation
\begin{equation}
\df \tau{t}=\Theta^{p-1}
\end{equation}
such that $\tau(0)=0$ we conclude that $\tau$ ranges from $0$ to $\un$ and $\tilde u$ satisfies
\begin{equation}
\pde{\tilde u}\tau=vF^{-p}-n^{-p}\tilde u,
\end{equation}
or equivalently, with $\tilde\f=\log\tilde u$
\begin{equation}\lae{5.7}
\pde{\tilde\f}{\tau}=v\tilde u^{-1}F^{-p}-n^{-p}.
\end{equation}
\el

Since the spatial derivatives of $\tilde \f$ are identical to those $\f$ equation \re{5.7} is a  nonlinear parabolic equation with a uniformly elliptic and concave operator $F$. Hence we can apply the Krylov-Safonov estimates and thereafter the parabolic Schauder estimates to conclude:
\bt
The rescaled flow
\begin{equation}
\tilde x'=F^{-p}\nu-n^{-p}\tilde x,
\end{equation}
where $\tilde x=x \Theta^{-1}$ and a prime indicates the  differentiation with respect to $\tau$, note that $\tilde x=\tilde x(\tau,\xi)$, exists for all time and the leaves converge in $C^{m+2,\bet}(\Ss[n])$, $0\le m\le \un$, $0<\bet<\al$, to a sphere of radius $1$ provided that $M_0\in C^{m+2,\al}$, $0<\al<1$,  is a star-shaped admissible initial hypersurface and the curvature function $F$ is of class $C^{m,\al}$, where we note that in case $1<p<\un$ we assume $\C=\C_+$ such that the admissible hypersurfaces are strictly convex.
\et


\bibliographystyle{hamsplain}

\begin{thebibliography}{10}

\bibitem{chow-gulliver}
Bennett Chow and Robert Gulliver, \emph{Aleksandrov reflection and nonlinear
  evolution equations. {I}. {T}he {$n$}-sphere and {$n$}-ball}, Calc. Var.
  Partial Differential Equations \textbf{4} (1996), no.~3, 249--264,
  {\href{http://dx.doi.org/10.1007/s005260050038}{doi:10.1007/s005260050038}}.

\bibitem{cg90}
Claus Gerhardt, \emph{Flow of nonconvex hypersurfaces into spheres}, J. Diff.
  Geom. \textbf{32} (1990), 299--314,
  {\href{http://math.uni-heidelberg.de/studinfo/gerhardt/jdg90.pdf}{pdf file}}.

\bibitem{cg:spaceform}
\bysame, \emph{Closed {Weingarten} hypersurfaces in space forms}, Geometric
  Analysis and the Calculus of Variations (J{\"urgen}~Jost, ed.), International
  Press, Boston, 1996,
  {\href{http://www.math.uni-heidelberg.de/studinfo/gerhardt/WeingartenInSpaceforms.pdf}{pdf
  file}}, pp.~71--98.

\bibitem{cg:cp}
\bysame, \emph{Curvature {P}roblems}, Series in Geometry and Topology, vol.~39,
  International Press, Somerville, MA, 2006.

\bibitem{cg:survey}
\bysame, \emph{Curvature flows in semi-{R}iemannian manifolds}, Geometric Flows
  (Huai-Dong Cao and Shing-Tung Yau, eds.), Surveys in Differential Geometry,
  vol. XII, International Press of Boston, Somerville, MA, 2008,
  {\href{http://arXiv.org/abs/0704.0236}{arXiv:0704.0236}}, pp.~113--165.

\bibitem{huisken:jdg1}
Gerhard Huisken, \emph{Flow by mean curvature of convex surfaces into spheres},
  J. Differential Geom. \textbf{20} (1984), no.~1, 237--266.

\bibitem{li:inverse}
Qi-Rui Li, \emph{Surfaces expanding by the power of the {G}auss curvature
  flow}, Proc. Amer. Math. Soc. \textbf{138} (2010), no.~11, 4089--4102,
  {10.1090/S0002-9939-2010-10431-8}.

\bibitem{mccoy:vol1}
James McCoy, \emph{The surface area preserving mean curvature flow}, Asian J.
  Math. \textbf{7} (2003), no.~1, 7--30.

\bibitem{oliver:2004surfaces}
O.C. Schn{\"u}rer, \emph{{Surfaces expanding by the inverse Gau{\ss} curvature
  flow}}, J. reine angew. Math. \textbf{600} (2006), 117--134,
  {\href{http://arXiv.org/abs/math/0412297}{arXiv:math/0412297}}.

\bibitem{urbas:convex}
John I.~E. Urbas, \emph{An expansion of convex hypersurfaces}, J. Differential
  Geom. \textbf{33} (1991), no.~1, 91--125.

\end{thebibliography}

\providecommand{\bysame}{\leavevmode\hbox to3em{\hrulefill}\thinspace}
\providecommand{\href}[2]{#2}



\end{document}